\documentclass[a4paper]{amsart}
\usepackage{amsmath,amsfonts,amsbsy,amsgen,amscd,mathrsfs,amssymb,amsthm,mathtools}
\usepackage[colorlinks=true,citecolor=blue,linkcolor=blue]{hyperref}

\usepackage{hhline}

\usepackage{accents}

\usepackage{pgfplots}
\usepackage{tikz-cd}
\usetikzlibrary{arrows,automata}

\numberwithin{equation}{section}
\newtheorem{theorem}{Theorem}[section]
\newtheorem{corollary}[theorem]{Corollary}

\newtheorem{lemma}[theorem]{Lemma}
\newtheorem{proposition}[theorem]{Proposition}
\theoremstyle{definition}

\newtheorem{definition}[theorem]{Definition}

\newtheorem{notation}[theorem]{Notation}
\newtheorem{problem}[theorem]{Problem}

\newtheorem{remark}[theorem]{Remark}

\makeatletter

\newcommand\be{\beta}
\newcommand\dd{\mathrm d}

\newcommand\de{\delta}
\newcommand\deq{\stackrel{\mathrm{distr.}}{=}}
\newcommand\eps{\varepsilon}

\newcommand\ka{\kappa}
\newcommand\La{\Lambda}

\newcommand\si{\sigma}

\newcommand\ze{\zeta}
\renewcommand\d{~\mathrm d}
\renewcommand\phi{\varphi}

\newcommand\mbb{\mathbb}
\newcommand\mbf{\mathbf}
\newcommand\mc{\mathcal}
\newcommand\mf{\mathfrak}
\newcommand\mr{\mathrm}

\begin{document}

\title{Number Rigidity of the Stochastic Airy Operator}
\author[P. Y. Gaudreau Lamarre]{Pierre Yves Gaudreau Lamarre}
\address{University of Chicago}
\email{pyjgl@uchicago.edu}
\author[P. Ghosal]{Promit Ghosal}
\address{Massachusetts Institute of Technology}
\email{promit@mit.edu}
\author[W. Li]{Wenxuan Li}
\address{University of Chicago}
\email{wenxuanl@uchicago.edu}
\author[Y. Liao]{Yuchen Liao}
\address{University of Warwick}
\email{Yuchen.Liao@warwick.ac.uk}
\subjclass[2010]{}

\begin{abstract}
We prove that the spectrum of the stochastic Airy operator is rigid in the sense of Ghosh and Peres \cite{GP} for Dirichlet and Robin boundary conditions. This proves the rigidity of the Airy-$\beta$ point process and the soft-edge limit of rank-$1$ perturbations of Gaussian $\beta$-Ensembles for any $\beta>0$, and solves an open problem mentioned in \cite{BNQ}. Our proof uses a combination of the semigroup theory of the stochastic Airy operator and the techniques for studying insertion and deletion tolerance of point processes developed in \cite{HS}. 
\end{abstract}

\maketitle

\section{Introduction}

In this paper, we are interested in the spectrum of a random Schr\"odinger operator
called the stochastic Airy operator. The stochastic Airy operator
acts on functions $f:[0,\infty)\to\mbb R$ on the positive half line
satisfying a general boundary condition of the form $af'(0)+bf(0)=0$ for some $a,b\in\mbb R$.
The operator in question, which we denote by $H_\be$, is defined as
\begin{align}
\label{Equation: SAO}
H_\be f(x):=-f''(x)+\Big(x+\tfrac2{\sqrt\be}W'(x)\Big)f(x),
\end{align}
where $\be>0$ is a fixed parameter and $W$ is a standard Brownian motion. Given that Brownian motion sample
paths are not differentiable, a rigorous definition of $H_\be$ is nontrivial. Nevertheless,
it can be shown that once a suitable boundary condition is fixed, $H_\be$ is a bonafide
semibounded self-adjoint operator with compact resolvent when defined through its quadratic form
\cite{BV,GL0,M,RRV}. We refer to Section \ref{Section: SAO Rigorous Definition} below for more details.

Much of the interest of studying the spectrum of $H_\be$ comes from the theories of
random matrices and interacting particle systems. Following the seminal work of Dumitriu
and Edelman \cite{DE}, Edelman and Sutton \cite{ES} introduced the stochastic Airy operator
as a means of characterizing the soft-edge scaling limits of interacting particle systems
known as the $\be$-Hermite and $\be$-Laguerre ensembles (see, e.g., \cite[Sections 2.1.1 and 2.1.2]{ES}).
This program was later realized by Ram\'irez, Rider
and Vir\'ag \cite{RRV}, and then extended to more general ensembles by Bloemendal,
Krishnapur, Rider, and Vir\'ag \cite{BV,KRV} (respectively, rank-1 spiked $\be$-ensembles
and general Dyson $\be$-ensembles).

In this paper, we prove that $H_\be$'s eigenvalue point process is number rigid
for any choice of $\be>0$, both for Dirichlet and Robin boundary conditions.
Roughly speaking, a point process $\Pi$ on $\mbb R^d$ is number rigid if
there exists a class of subsets $B\subset\mbb R^d$ such that
the number of $\Pi$'s points inside $B$ is uniquely determined by the configuration
of points outside $B$. We prove that this rigidity property holds for
$H_\be$'s spectrum for every Borel set $B\subset\mbb R$ that is bounded
above, that is, for which there exists some $c\in\mbb R$ such that
$B\subset(-\infty,c]$. We refer to Theorem \ref{Theorem: Main} below for a precise statement.

The result proved in this paper contributes to the program that is concerned with understanding the
spatial conditioning of point processes with dependence. That is, if we condition a point process
on having some fixed configuration outside of some set $B$, then what is the conditional distribution
of points inside $B$? In the specific context of random matrix theory, such problems have been
motivated by imperfect observations in nuclear physics experiments and the connection between
energy levels of heavy nuclei and random matrix eigenvalues; see \cite[Section 1.1]{CG} and
references therein for more details. From this point of view, number rigidity can be seen as a first
step towards exploring questions of spatial conditioning for edge $\be$-Ensembles, and naturally
leads to further investigations in this direction, some of which we discuss in Section \ref{Section: Open Prob.}.

\subsection{Previous Results}

The concept of number rigidity was formally introduced by Ghosh and Peres in \cite{GP}
(see also \cite{AM,HS}) in the context of the Ginibre process and Gaussian analytic function
zeroes, and has since then led to a number of interesting developments. Using a variety of
methods, number rigidity has been established for a number of point
processes, including for example perturbed lattices \cite{HS,PS}, determinantal/Pfaffian
point processes \cite{B,BNQ}, hyperuniform point processes \cite{GL},
Sine$_\be$ point processes \cite{DHLM,G,CN}, and more.

Looking more specifically at the stochastic Airy operator, until now the number rigidity
of $H_\be$'s spectrum was only known in the case $\be=2$ with Dirichlet boundary
condition (i.e., $f(0)=0$). This result, which is due to Bufetov \cite{B}, relies crucially on the fact that
in this particular case, $H_\be$'s eigenvalues form a determinantal point process. Thus,
the method used in \cite{B} is not amenable to extension to $\be\neq2$. The number rigidity of $H_\be$'s spectrum
with $\be=1,4$ and Dirichlet boundary condition was raised as an open problem in \cite{BNQ}
in the context of Pfaffian point processes; to the best of our knowledge,
this question remained open until now.

In particular, the methods previously used to prove
number rigidity for general (i.e., $\be\neq2$) Sine-$\be$ processes in \cite{DHLM,CN} do not appear to be
directly applicable to the edge.
On the one hand, the method used in \cite{CN} relies on moment estimates of the circular $\beta$-Ensembles,
from which the Sine-$\be$ process can be obtained as a scaling limit.
On the other hand, the proof in \cite{DHLM} relies on establishing a much stronger property
encapsulated by the canonical DLR equation (\cite[Theorem 2.1]{DHLM}). The proof of this DLR property in \cite{DHLM} also relies
on using the circular $\be$-Ensemble approximation, which is not available at the edge.
Rigidity properties of general $\beta$-Ensembles is addressed in \cite{BEY14} and further refined in \cite{AH20,BMP21}; however, these estimates of the rigidity do not pass to the limit under the edge scaling and thus do not appear
to immediately yield number rigidity for the Airy-$\beta$ process.

Motivated in part by the desire to prove the number rigidity of $H_\be$ for $\be\neq2$,
in \cite{GLGL} a new strategy to prove number
rigidity in the spectrum of general random Schr\"odinger operators was developed.
The strategy in question, which is based on a sufficient condition for rigidity
due to Ghosh and Peres \cite{GP}, consists of proving that the variance of the trace
of the semigroup $\mr e^{-t H}$ of a Schr\"odinger operator $H$ vanishes as $t\to0$
using the Feynman-Kac formula. We refer to Section \ref{Section: Semigroup Approach}
for more details. While the main result of \cite{GLGL} establishes the number
rigidity of the spectrum of a wide class of one-dimensional continuous random Schr\"odinger operators,
it unfortunately does not apply to $H_\be$, as in this case the variance of the
trace of $H_\be$'s semigroup is merely bounded as $t\to0$.

In this context, the main improvement in this paper consists of replacing the sufficient
condition for number rigidity by Ghosh and Peres with a more general result
inspired by the work of Holroyd and Soo \cite{HS}. The result in question,
which we state as Lemma \ref{Lemma: New Condition} in the context
of the stochastic Airy operator, allows to replace the vanishing of the trace of
the semigroup by a weaker covariance condition.

Lastly, we note that the rigidity properties of point processes is
an important tool in the derivation of a variety of fundamental results in random matrix theory
and interacting particle systems,
such as local laws and universality.
The most commonly-used notion of rigidity considers the maximal deviation of the point process from a set of deterministic
values that, in some sense, captures the typical behavior of each point.
Many of the recent works in this direction have focused on the rigidity
of general $\beta$-Ensembles (e.g, \cite{AH20,BEY14,BMP21,CFLW21}).
The results of this paper studies the rigidity of $\be$-Ensembles from a different point of view---the number rigidity of the edge-limit---and thus complements these previous investigations (in particular, number rigidity is not a corollary of this other notion of
rigidity).

\subsection{Open Problems}
\label{Section: Open Prob.}

The results in this paper raise a number of interesting questions for future research.
Most notably, given that number rigidity is now established for both the bulk and soft-edge limits of $\be$-Ensembles
and their rank-1 perturbations (by combining \cite{B,BNQ,DHLM,G,CN} with the present paper), it is natural to investigate
whether or not number rigidity occurs in the other point processes that describe the scaling limits of random matrix ensembles:

\begin{problem}
For a given integer $r\geq2$ and a rank-$r$ spike matrix
$S=\sum_{i=1}^rw_iu_iu_i^\dagger$,
consider the $r\times r$ multivariate stochastic Airy operator
\[H^{(r)}_\be:=-\tfrac{\dd^2}{\dd x^2}+\sqrt 2B'_x+rx\]
with initial condition $f'(0)=Sf(0)$,
as defined in \cite{BV2} for $\be=1,2,4$ (i.e., the soft-edge scaling limit of rank-$r$
perturbations of Gaussian $\be$-Ensembles). Is the eigenvalue point
process of $H^{(r)}_\be$ number rigid?
\end{problem}

\begin{problem}
Consider the hard-edge scaling limit of the $\be$-Laguerre ensembles,
as characterized by the stochastic Bessel operator or its spiked version
(i.e., \cite{ES,RR1,RR2,RR3}). Is this point process number rigid?
\end{problem}

While it seems plausible that these problems could be approached with semigroup methods, substantial
modifications would need to be implemented in order to extend our arguments to this setting. We thus leave
these two problems open for future investigations.

Finally, it would be interesting to see if more detailed information
about the spatial conditioning of the general Airy-$\be$ process could be obtained. In particular, a natural starting point would be the following:

\begin{problem}
\label{Problem: Spatial Conditioning}
Fix a boundary condition for $H_\be$, and let $B\subset\mbb R$
be a Borel set that is bounded above.
Suppose that we condition on the configuration of $H_\be$'s eigenvalues
outside $B$. By number rigidity, there exists a deterministic $N\in\mbb N$
(which only depends on the configuration of eigenvalues outside $B$)
that is equal to the number of $H_\be$'s eigenvalues inside $B$.
Let $\ze\in B^N$ be the random vector whose components are $H_\be$'s eigenvalues inside $B$
(conditional on the configuration outside $B$), taken in a uniformly random order.
What is the support of $\ze$'s probability distribution on $B^N$?
In particular, is $H_\be$'s eigenvalue point process {\it tolerant} in the sense of Ghosh and Peres \cite{GS}; that is,
are $\ze$'s probability distribution and the Lebesgue measure on $B^N$ mutually absolutely continuous?
\end{problem}

Number rigid point processes for which tolerance has been proved include the
Ginibre point process \cite[Theorem 1.2]{GS} and the Sine-$\be$ process for every $\be>0$ \cite[Corollary 1.2]{DHLM}
(see also \cite{BR} for the Sine-$2$ process).
To the best of our knowledge, no such result is known for the Airy-$\be$ process for any choice of $\be>0$.

\subsection{Organization}

The remainder of this paper is organized as follows:
In Section \ref{Section: Main Result}, we provide a precise
definition of the stochastic Airy operator and state our
main result, namely, Theorem \ref{Theorem: Main}. In Section \ref{Section: Outline}, we provide
an outline of the main steps of our proof of Theorem \ref{Theorem: Main}.
Then, in Section \ref{Section: Technical}, we provide the proofs
of a variety of technical results stated in Section \ref{Section: Outline}.

\section{Setup and Main Result}
\label{Section: Main Result}

\subsection{Stochastic Airy Operator}
\label{Section: SAO Rigorous Definition}

The main goal of this section is to introduce the stochastic Airy operator and state the main result of this paper. We begin by introducing a few notations.   
We use $V:[0,\infty)\to\mbb R$ to denote the linear function $V(x)=x$. Given $f,g\in L^2\big([0,\infty)\big)$, we use the standard notations
\[\langle f,g\rangle:=\int_0^\infty f(x)g(x)\d x
\qquad\text{and}\qquad
\|f\|_2:=\left(\int_0^\infty f(x)^2\d x\right)^{1/2}.\]
We define the norm
\[\|f\|_*^2:=\|f'\|_2^2+\|V^{1/2}f\|_2^2+\|f\|_2^2.\]
We use $C^\infty_0\big([0,\infty)\big)$ to denote the set of smooth and compactly
supported functions from $[0,\infty)$ to $\mbb R$.

Let $W$ be a standard Brownian motion.
Given a continuous and compactly supported map $\phi:[0,\infty)\to\mbb R$, we denote
the stochastic integral
\[\xi_\be(\phi):=\frac2{\sqrt\be}\int_0^\infty \phi(x)\d W(x)\]
for every $\be>0$.
In particular, $\phi\mapsto\xi_\be(\phi)$ is Gaussian process with covariance
given by the inner product $\frac4\be\langle\cdot,\cdot\rangle$.
The basis of the definition of $H_\be$'s spectrum is the following quadratic form:

\begin{definition}
Let $\be>0$ and $w\in(-\infty,\infty]$.
Define the set
\[\mr{FC}_{\be,w}:=\begin{cases}
C^\infty_0\big([0,\infty)\big)&\text{if }w<\infty\\
\big\{\phi\in C^\infty_0\big([0,\infty)\big):\phi(0)=0\big\}&\text{if }w=\infty.
\end{cases}\]
Then, for every $\phi\in\mr{FC}_{\be,w}$, we define the quadratic form
\begin{align}
\label{Equation: SAO Form}
\mc E_{\be,w}(\phi):=
\begin{cases}
-w\,\phi(0)^2+\|\phi'\|_2^2+\|V^{1/2}\phi\|_2^2+\xi_\be(\phi^2)&\text{if }w<\infty\\
\|\phi'\|_2^2+\|V^{1/2}\phi\|_2^2+\xi_\be(\phi^2)&\text{if }w=\infty.
\end{cases}
\end{align}
\end{definition}

The spectrum of the stochastic Airy operator was first constructed using the form \eqref{Equation: SAO Form}
in \cite{BV,RRV}. The following statement is a special case of \cite[Propositions 3.2 and 3.4]{GL0}
(see also \cite{M}):

\begin{theorem}[\cite{BV,GL0,M,RRV}]
Let $\be>0$ and $w\in(-\infty,\infty]$.
The quadratic form $\mc E_{\be,w}$ can be continuously extended with respect
to $\|\cdot\|_*$ to a domain $D(\mc E_{\be,w})$ on which
it is closed and semibounded. In particular, there exists a unique self-adjoint
operator $H_{\be,w}$ on $L^2\big([0,\infty)\big)$ whose quadratic form
is given by \eqref{Equation: SAO Form}. Almost surely, $H_{\be,w}$ has compact
resolvent.
\end{theorem}

By examining \eqref{Equation: SAO Form}, we note that $w<\infty$
corresponds to the case where $H_\be$ acts on functions with Robin
boundary condition $f'(0)=wf(0)$, and $w=\infty$ corresponds to the
Dirichlet boundary condition $f(0)=0$.
As an immediate consequence of the above result, we conclude the following regarding
$H_{\be,w}$'s spectrum:

\begin{corollary}
\label{Corollary: SAO}
Let $\be>0$ and $w\in(-\infty,\infty]$. Almost surely, the operator $H_{\be,w}$
has a purely discrete spectrum of eigenvalues
\begin{align}
\label{Equation: SAO Spectrum}
-\infty<\La_{\be,w}(1)\leq\La_{\be,w}(2)\leq\La_{\be,w}(3)\leq\cdots
\end{align}
that is bounded below and without accumulation point.
\end{corollary} 

\begin{remark}
The eigenvalues can be defined through the
min-max formula using the quadratic form \eqref{Equation: SAO Form}. Following \cite{BV,RRV}, $\La_{\be,w}(1)$ can be described as
\begin{align}
\La_{\be,w}(1) := \inf\big\{\langle f, H_{\be,w} f\rangle: f\in D(\mc E_{\be,w}), \|f\|_2= 1\big\}.
\end{align}
Suppose $f^{*}_1\in \mc E_{\be,w}$ is such that $\La_{\be,w}(1)$ is equal to $\langle f_1^*, H_{\be,w} f_1^*\rangle$. Then, $\La_{\be,w}(2)$ is the infimum of the quadratic form over $f\in D(\mc E_{\be,w})$ with $\|f\|_2=1$ that are orthogonal to $f^{*}_1$ w.r.t. $\langle \cdot, \cdot \rangle$, and so on for higher eigenvalues.
\end{remark}

\subsection{Main Result}

We first define the eigenvalue point process of $H_{\beta, w}$
and then state our main result.

\begin{definition}[Point Process]
Let $A\subset\mbb R$ be a Borel set
and $\mc B_A$ be the associated Borel $\si$-algebra.
We use $\mc N^\#_A$ to denote the set
of integer-valued measures $\mu$ on $(A,\mc B_A)$ such that $\mu(K)<\infty$ for every bounded set $K\subset A$.
We equip $\mc N^\#_A$ with the topology generated by the maps
\[\mu\mapsto\int f(x)\d\mu(x),\qquad f\text{ continuous and compactly supported,}\]
as well as the associated Borel $\si$-algebra.
For every $\mu\in\mc N^\#_A$, there exists a sequence
$(x(k))_{k\in\mbb N}\subset A$ without accumulation point such that
$\mu=\sum_{k=1}^\infty\de_{x(k)}$,
where $\de_x$ denotes the Dirac mass at a point $x\in\mbb R$;
that is, for any Borel set $B\subset\mbb R$, one has
$\de_x(B)=\mbf 1_{\{x\in B\}}$.
A point process on $A$ is a random element
that takes values in $\mc N^\#_A$.
\end{definition}

\begin{definition}[Eigenvalue point process and its restrictions]
For every $\be>0$ and $w\in(-\infty,\infty]$, we define
$H_{\be,w}$'s eigenvalue point process as
\[\Pi_{\be,w}:=\sum_{k=1}^\infty\de_{\La_{\be,w}(k)}\in\mc N^\#_\mbb R.\]
For every Borel set $A\subset\mbb R$, we denote the restriction
\[\Pi_{\be,w}|_{A^c}:=\sum_{k=1}^\infty\mbf 1_{\{\La_{\be,w}(k)\in A^c\}}\,\de_{\La_{\be,w}(k)}\in\mc N^\#_{A^c}.\]
That is, the same random measure as $\Pi_{\be,w}$, but excluding the masses in $A$.
\end{definition}

We say that a Borel set $B\subset\mbb R$ is \emph{bounded above} if
there exists some $c\in\mbb R$ such that $B\subset(-\infty,c]$.
By Corollary \ref{Corollary: SAO},
$\Pi_{\be,w}$ is a random counting measure on $\mbb R$
such that, almost surely, $\Pi_{\be,w}(B)<\infty$ for every
Borel set $B$ that is bounded above. With this in hand, the main result of this paper is the following:

\begin{theorem}
\label{Theorem: Main}
For every $\be>0$, $w\in(-\infty,\infty]$, and any bounded above Borel set $B\subset\mbb R$, there exists
a deterministic measurable function $F_B:\mc N^\#_{B^c}\to\mbb R$ such that
\begin{align}\label{eq:PiRigidity}
\Pi_{\beta, w}(B)=F_B\big(\Pi_{\beta, w}|_{B^c}\big)
\end{align}
almost surely, where $\Pi_{\beta,w}(B)$ is the number of points of $\Pi_{\beta, w}$ inside the Borel set $B$ and $\Pi_{\be,w}|_{B^c}\in\mc N^\#_{B^c}$ is the restriction of the point process $\Pi$ to the complement of $B$.
\end{theorem}

\begin{remark}
Theorem~\ref{Theorem: Main} shows that the number of $H_{\beta, w}$'s eigenvalues
inside $B$ is uniquely determined by the configuration of $H_{\beta,w}$'s eigenvalues
outside of $B$. As per the terminology introduced by Ghosh and Peres in \cite{GP}, any point process which satisfies this property for any bounded Borel set is number rigid. Since any bonded Borel set on $\mathbb{R}$ is bounded above, Theorem~\ref{Theorem: Main} proves that $H_{\beta, w}$ is number rigid in the sense of Ghosh and Peres \cite{GP}.    
\end{remark}

\section{Proof Outline}
\label{Section: Outline}

In this section, we provide an outline of the proof of Theorem \ref{Theorem: Main}.
In order to keep the argument readable, we defer the proof of several technical
lemmas to Section \ref{Section: Technical}. We also take this opportunity to showcase
how the proof method used in this paper improves on that of \cite{GLGL}.

\begin{remark}
Unless otherwise mentioned, for the remainder of this section we assume that
$\be>0$ and $w\in(-\infty,\infty]$ are fixed and arbitrary.
\end{remark}

\subsection{Step 1 - The Semigroup Approach: A New Sufficient Condition}
\label{Section: Semigroup Approach}

Our method to prove Theorem~\ref{Theorem: Main} is based on the semigroup of the stochastic Airy operator. More specifically, the
key objects of study are the traces
\[\mr{Tr}[\mr e^{-t H_{\be,w}}]=\sum_{k=1}^\infty\mr e^{-t \La_{\be,w}(k)}=\int_{\mbb R}\mr e^{-t x}\d\Pi_{\be,w}(x)\]
for some positive $t>0$. Our proof of rigidity is based on the following sufficient condition:

\begin{lemma}
\label{Lemma: New Condition}
Suppose that
\[\limsup_{t\to0}\mbf{Var}\big[\mr{Tr}[\mr e^{-t H_{\be,w}}]\big]<\infty,\]
and that for every fixed $u>0$,
\begin{align}
\label{Equation: Vanishing Covariance}
\lim_{t\to0}\mbf{Cov}\big[\mr{Tr}[\mr e^{-tH_{\be,w}}],\mr{Tr}[\mr e^{-uH_{\be,w}}]\big]=0.
\end{align}
Then, $H_{\be,w}$ satisfies \eqref{eq:PiRigidity} for any bounded above Borel set $B$. 
\end{lemma}

Lemma \ref{Lemma: New Condition} is proved in Section \ref{Section: Sufficient Condition}.
This result is inspired by the work of Holroyd and Soo \cite{HS} on {\it deletion-tolerance}
(see \cite[Page 2]{HS} for a definition of deletion-tolerance), as well as the connection between deletion-tolerance and number rigidity
that was pointed out in \cite[Proposition 1.2]{PS}.
More specifically, a sufficient condition similar to Lemma \ref{Lemma: New Condition}
was used in \cite[Lemma 7.2]{HS} to prove that the two-dimensional perturbed lattice
is not deletion-tolerant.
We note that in \cite[Lemma 7.2]{HS}, the trace $\mr{Tr}[\mr e^{-t H_{\be,w}}]$
is replaced by a different linear functional of the perturbed lattice,
denoted $\La(h_k)$ therein; see Remark \ref{Remark: Linear Statistics} for more details
on the distinction between traces of semigroups and general linear statistics
in questions of rigidity.
We refer to Section \ref{Section: Sufficient Condition} for the details of how the argument used in
\cite[Lemma 7.2]{HS} is adapted to our setting.

Semigroup theory is an attractive strategy to prove number rigidity
for general Schr\"odinger operators (including $H_{\be,w}$) because of the Feynman-Kac formula.
In short, the Feynman-Kac formula states that the semigroups of Schr\"odinger
operators admit an explicit formulation in terms of elementary stochastic processes; see, e.g.,
\cite{S}. The Feynman-Kac formula for $H_{\be,w}$ specifically was proved in \cite{GS} in the
case $w=\infty$ and \cite{GLS} in the case $w<\infty$.  Following \cite[Theorem 2.24]{GL0}, we first introduce the necessary notations, and then proceed to state the Feynman-Kac formula for the trace of $\mr e^{-tH_{\be,w}}$.

Let $B$ denote a standard Brownian motion on $\mbb R$, and
Let $X$ denote a reflected Brownian motion on $(0,\infty)$
(i.e., the absolute value of a Brownian motion).
Let $Z=B$ or $Z=X$;
for every $x,y,t>0$ we denote the conditioned process
\[Z^x:=\big(Z\big|Z(0)=x\big)
\qquad\text{and}\qquad
Z^{x,y}_t:=\big(Z\big|Z(0)=x\text{ and }Z(t)=y\big).\]

We use $L_t^a(Z)$ ($a,t>0$) to denote the continuous version of the local time of some process $Z$ (or its conditioned versions) on the time interval $[0,t]$, that is, the stochastic process such
that for every Borel measurable $f:(0,\infty)\to\mbb R$, one has
\[\int_0^tf\big(Z(s)\big)\d s=\int_0^\infty L_t^a(Z)f(a)\d a=\langle L_t(Z),f\rangle.\]
(See, e.g., \cite[Chapter VI, Corollary 1.6 and Theorem 1.7]{RY}.)

We use $\mathfrak L_t^0(Z)$ ($t>0$) to denote the boundary local time
of a process $Z$ (or its conditioned versions) on the time interval $[0,t]$, that is,
\[\mathfrak L_t^0(Z):=\lim_{\eps\to0}\frac1{2\eps}\int_0^t\mbf 1_{\{Z(s)\in[0,\eps)\}}\d s.\]

\begin{theorem}[\cite{GL0,GLS,GS}]\label{Theorem: Trace Identity}
For every $t>0$,
\begin{align}
\label{Equation: Trace Identity}
\mr{Tr}[\mr e^{-t H_{\be,w}/2}]
=\int_0^\infty\frac{1+\mr e^{-2x^2/t}}{\sqrt{2\pi t}}\mbf E\left[\mr e^{-\langle L_t(X^{x,x}_t),\frac12V\rangle-\frac12\xi_\be(L_t(X^{x,x}_t))-w\mf L_t^0(X^{x,x}_t)}\Big|\xi_\be\right]\d x,
\end{align}
where $X$ is assumed to be independent of $\xi_\be$ in the above conditional expectation.
We use the convention $\mr e^{-\infty\cdot\mf L_t^0(X)}=\mbf 1_{\{\mf L_t^0(X)=0\}}$
in the case where $w=\infty$, and we recall that $V(x)=x$.
\end{theorem}

In Sections \ref{Section: Step 2 Outline} and \ref{Section: Step 3 Outline}, we outline how
Theorem \ref{Theorem: Trace Identity} can be used to show that the conditions in
Lemma \ref{Lemma: New Condition} hold, and thus that $H_{\be,w}$'s spectrum is number
rigid. However, before moving on to this outline, we briefly comment on how Lemma \ref{Lemma: New Condition}
improves on our previous work in \cite{GLGL}.
In \cite{GLGL}, the sufficient condition for number rigidity that we used
was the following special case of the general strategy devised by Ghosh and Peres in
\cite[Theorem 6.1]{GP}:

\begin{proposition}[\cite{GP}]
\label{Proposition: Old Condition}
Let $H$ be a random Schr\"odinger operator
with (almost surely) compact resolvent.
If
\begin{align}\label{eq:Final Limit}
\lim_{t\to0}\mbf{Var}\big[\mr{Tr}[\mr e^{-t H}]\big]=0,
\end{align}
then $H$'s spectrum is number rigid.
\end{proposition}

\begin{remark}
\label{Remark: Linear Statistics}
The strategy in \cite[Theorem 6.1]{GP} is more general than what is stated in Proposition \ref{Proposition: Old Condition}:
Let $\Pi$ be a point process.
To prove that $\Pi$ is number rigid, it suffices to show that for any bounded Borel set $B\subset\mbb R$, there is a sequence of functions $(f_n)_{n\in\mbb N}$ that converge
uniformly to 1 on $B$ such that
the variance of the {\it linear statistics} $\mbf{Var}\big[\int f_n(x)\d\Pi(x)\big]$ vanishes as $n\to\infty$.
In this paper and \cite{GLGL}, we use the functions $f(x)=\mr e^{-t x}$ ($t>0$) because the Feynman-Kac
formula makes it possible to compute explicitly the variance and covariance of the trace of $\mr e^{-t H}$.
\end{remark}

Unfortunately, the best that we could achieve in \cite{GLGL} using \eqref{Equation: Trace Identity} was the following result:

\begin{proposition}[\cite{GLGL}]
\label{Proposition: Bounded Variance}
It holds that
\begin{align}
\label{Equation: Bounded Variance}
\limsup_{t\to0}\mbf{Var}\big[\mr{Tr}[\mr e^{-t H_{\be,w}}]\big]<\infty.
\end{align}
\end{proposition}

Using the notation/terminology of \cite{GLGL}, the above proposition follows from
\cite[Theorem 4.1]{GLGL} in {\bf Case 2} (i.e., operators on the positive half-line)
with white noise (whereby $\mf d=3/2$ by \cite[(2.16)]{GLGL}),
and $\mf a=1$, $\ka=1/2$, and $\nu=0$, the latter of which corresponds to the case $\frac12 V=\frac x2$.
Moreover, by exploiting the special integrable structure available in the case $\be=2$, it can be
shown that the finiteness of the limsup in \eqref{Equation: Bounded Variance} cannot be improved
to a vanishing limit in general:

\begin{proposition}[{\cite[Proposition 2.27]{GLGL}}]
\label{prop:GLGL}
It holds that
\[\lim_{t\to0}\mbf{Var}\big[\mr{Tr}[\mr e^{-t H_{2,\infty}}]\big]=(4\pi)^{-1}.\]
\end{proposition}

Proposition~\ref{prop:GLGL} shows that \eqref{eq:Final Limit} cannot hold for $H_{\be,w}$ when $\be=2$ and $w= \infty$.
Thus, the general number rigidity of $H_{\be,w}$ cannot be proved using  Proposition \ref{Proposition: Old Condition} since the condition \eqref{eq:Final Limit} cannot be met for every choice of $\be$ and $w$. 
In particular, one of the key innovations in this paper is to replace the vanishing variance condition \eqref{eq:Final Limit}
with the weaker conditions \eqref{Equation: Vanishing Covariance} and \eqref{Equation: Bounded Variance}.
We now explain how we use the Feynman-Kac formula to achieve these weaker conditions.

\subsection{Step 2 - Covariance Formula}
\label{Section: Step 2 Outline}

Thanks to \eqref{Equation: Bounded Variance} and Lemma \ref{Lemma: New Condition},
the proof of Theorem \ref{Theorem: Main} is reduced to proving \eqref{Equation: Vanishing Covariance}.
For this purpose, we derive a formula of the covariance in \eqref{Equation: Vanishing Covariance}. Before proceeding to that derivation, we introduce a few notations.

For every $t,u,x,y>0$, we let $X^{x,x}_t$ and $\bar X^{y,y}_u$
be independent, where $\bar X^{y,y}_u$ has the same law as $X^{y,y}_u$.
We define the shorthands
\begin{align*}
\mc P_{t,u}(x,y)&:=\left(\frac{1+\mr e^{-2x^2/t}}{\sqrt{2\pi t}}\right)\left(\frac{1+\mr e^{-2y^2/u}}{\sqrt{2\pi u}}\right), \\
\mc A_{t,u}(x,y) &:=-\langle L_t(X^{x,x}_t)+L_u(\bar X^{y,y}_u),\tfrac12V\rangle\\
\mc B_{t,u}(x,y) &:=-w\mf L_t^0(X^{x,x}_t)-w\mf L_u^0(\bar X^{y,y}_u), \\
\mc C_{t,u}(x,y) &:=\tfrac1{2\be}\big(\|L_t(X^{x,x}_t)\|_2^2+\|L_u(\bar X^{y,y}_u)\|_2^2\big)\\
\mc D_{t,u}(x,y) &:= \frac{\langle L_t(X^{x,x}_t),L_u(\bar X^{y,y}_u)\rangle}{\be}.
\end{align*}

In the following lemma, we state the formula of the covariance in \eqref{Equation: Vanishing Covariance} in terms of the above notations. 

\begin{lemma}
\label{Lemma: Covariance Formula}
For every $t,u>0$,
\[\mbf{Cov}\big[\mr{Tr}[\mr e^{-tH_{\be,w}/2}],\mr{Tr}[\mr e^{-uH_{\be,w}/2}]\big]
=\int_{(0,\infty)^2}\mc P_{t,u}(x,y)\mbf E\Big[\mr e^{\mc A_{t,u}(x,y)+\mc B_{t,u}(x,y)+\mc C_{t,u}(x,y)}\big(\mr e^{\mc D_{t,u}(x,y)}-1\big)\Big]\d x\dd y.\]
\end{lemma}

Lemma \ref{Lemma: Covariance Formula} is proved in Section \ref{Section: Covariance Formula}.
Applying H\"older's inequality on the right hand side of the above identity, we thus obtain that
\[\mbf{Cov}\big[\mr{Tr}[\mr e^{-tH_{\be,w}/2}],\mr{Tr}[\mr e^{-uH_{\be,w}/2}]\big]\]
is bounded above by
\begin{multline}
\label{Equation: Covariance Holder}
\sup_{x,y>0}\left\{\mc P_{t,u}(x,y)
\mbf E\left[\mr e^{4\mc B_{t,u}(x,y)}\right]^{1/4}
\mbf E\left[\mr e^{4\mc C_{t,u}(x,y)}\right]^{1/4}\right\}\\
\cdot\int_{(0,\infty)^2}\mbf E\left[\mr e^{4\mc A_{t,u}(x,y)}\right]^{1/4}
\mbf E\left[\left(\mr e^{\mc D_{t,u}(x,y)}-1\right)^4\right]^{1/4}
\d x\dd y.
\end{multline}

The following step, which is the last step of the proof, focuses on showing that the expression in \eqref{Equation: Covariance Holder} vanishes as $t\to0$. 

\subsection{Step 3 - Technical Estimates and Conclusion}
\label{Section: Step 3 Outline}

In order to provide a quantitative control on the terms appearing in
\eqref{Equation: Covariance Holder}, we need a number of technical estimates.
First, it is easy to see by definition of $\mc P$ that for every fixed $u>0$,
\begin{align}
\label{Equation: Transition Kernel Estimate}
\limsup_{t\to0}t^{1/2}\sup_{x,y>0}\mc P_{t,u}(x,y)<\infty.
\end{align}
For the remaining terms in \eqref{Equation: Covariance Holder},
we have the following:

\begin{lemma}
\label{Lemma: Self-Intersection}
For every $u>0$, there exists constant $c>0$ such that
\[\limsup_{t\to0}\sup_{x,y>0}\mbf E\left[\mr e^{4\mc B_{t,u}(x,y)}\right]\leq c<\infty\]
and
\[\limsup_{t\to0}\sup_{x,y>0}\mbf E\left[\mr e^{4\mc C_{t,u}(x,y)}\right]\leq c<\infty.\]
\end{lemma}

\begin{lemma}
\label{Lemma: Gaussian Term}
For every $u>0$, there exists a constant $c>0$ such that for every $t\in(0,1]$ and $x,y>0$, one has
\[\mbf E\left[\left(\mr e^{\mc D_{t,u}(x,y)}-1\right)^4\right]^{1/4}\leq c\mr e^{-(x-y)^2/c}\mbf E\left[\left(\mr e^{\mc D_{t,u}(x,y)}-1\right)^8\right]^{1/8}.\]
\end{lemma}

\begin{lemma}
\label{Lemma: Potential Term}
For every $u,c>0$,
\[\limsup_{t\to0}\int_{(0,\infty)^2}\mbf E\left[\mr e^{4\mc A_{t,u}(x,y)}\right]^{1/4}\mr e^{-(x-y)^2/c}\d x\dd y<\infty.\]
\end{lemma}

\begin{lemma}
\label{Lemma: Vanishing Term}
For every $u>0$,
\[\limsup_{t\to0}t^{-3/4}\sup_{x,y>0}\mbf E\left[\left(\mr e^{\mc D_{t,u}(x,y)}-1\right)^8\right]^{1/8}<\infty.\]
\end{lemma}

At this point, we insert \eqref{Equation: Transition Kernel Estimate}
and Lemma \ref{Lemma: Self-Intersection} in the first line of \eqref{Equation: Covariance Holder},
and then apply Lemmas \ref{Lemma: Gaussian Term}--\ref{Lemma: Vanishing Term}
to the integral on the second line on \eqref{Equation: Covariance Holder}.
This yields that for every $u>0$,
\[\mbf{Cov}\big[\mr{Tr}[\mr e^{-tH_{\be,w}}],\mr{Tr}[\mr e^{-uH_{\be,w}}]\big]=O_u(t^{-1/2}\cdot t^{3/4})=O_u(t^{1/4}),\]
as $t\to0$, where the constant in $O_u$ depends on $u$. This then proves Theorem \ref{Theorem: Main}.

\section{Technical Results}
\label{Section: Technical}

We now conclude the proof of Theorem \ref{Theorem: Main} by providing proofs
for the technical estimates stated in the previous section. Many of these results
follow from the analysis performed in \cite{GLGL}; we include detailed
references when such is the case.

\subsection{Proof of Lemma \ref{Lemma: New Condition}}
\label{Section: Sufficient Condition}

This proof draws inspiration from \cite[Proposition 7.1 and Lemma 7.2]{HS}, which the authors had used to study the insertion and deletion tolerance of point processes. We nevertheless provide the argument in full for the convenience of the readers.
The argument in question is based on the following generalization of the $L^2$ weak law 
of large numbers:

\begin{lemma}
\label{Lemma: WLLN}
  Let $X_1,X_2,\ldots$ be a sequence of random variables with $\mbf{E}[X_i] = 0$. Suppose there is 
  $C\in(0,\infty)$ so that $\sup_{i\in\mathbb{N}}\mbf{E}[X_i^2]\leq C$. Moreover, assume that 
  for each fixed $i$, we have $\mbf{E}[X_iX_j]\rightarrow 0$ as $j\rightarrow \infty$. Then there exists a subsequence 
  $(r_n)_{n\in\mathbb{N}}$ so that $(X_{r_1}+\cdots+X_{r_n})/n$ converges to zero in probability as $n\to\infty$. 
\end{lemma}

We provide a proof of Lemma \ref{Lemma: WLLN} in Section \ref{Section: WLLN}.
Before doing so, however, we explain how Lemma \ref{Lemma: WLLN} is used to complete the proof of Lemma~\ref{Lemma: New Condition}.

\subsubsection{Proof of Lemma~\ref{Lemma: New Condition}}

To prove Lemma \ref{Lemma: New Condition}, we first fix a bounded above Borel set $B\subset \mathbb{R}$.
Note that for any $t>0$, we can write
\begin{align}
\label{Equation: Split the Integral}
\mr{Tr}\left[\mr e^{-tH_{\beta,w}}\right]=\int_\mbb R\mr e^{-tx}\d\Pi_{\beta,w}(x)=\int_B \mr e^{-tx}\d\Pi_{\beta,w}(x)+\int_{B^c}\mr e^{-tx}\d\Pi_{\beta,w}(x).
\end{align}
Let $t_n$ be a sequence that converges to zero. 
By adding and subtracting $\frac{1}{n}\sum_{i=1}^n\mr{Tr}\left[\mr e^{-t_i H_{\beta,w}}\right]$ and
$\frac{1}{n}\sum_{i=1}^n\mbf{E}\left[\mr{Tr}[\mr e^{-t_iH_{\beta,w}}]\right]$,
and then split the traces $-\mr{Tr}\left[\mr e^{-t_i H_{\beta,w}}\right]$ using \eqref{Equation: Split the Integral}, we can write
$$\Pi_{\beta,w}(B) = \int_B 1\d\Pi_{\beta,w}(x)$$
as the sum of the following three terms:
\begin{align}
  \label{Equation: First Term of Inside} 
  \int_B 1\d\Pi_{\beta,w}(x)-\frac{1}{n}\sum_{i=1}^n\int_B \mr e^{-t_ix}\d\Pi_{\beta,w}(x)\\
  \label{Equation: Second Term of Inside}
  \frac{1}{n}\sum_{i=1}^n\left(\mr{Tr}\left[\mr e^{-t_i H_{\beta,w}}\right]-\mbf{E}\left[\mr{Tr}[\mr e^{-t_iH_{\beta,w}}]\right]\right)\\
  \label{Equation: Third Term of Inside}
  \frac{1}{n}\sum_{i=1}^n\left(\mbf{E}[\mr{Tr}[\mr e^{-t_iH_{\beta,w}}]]-\int_{B^c}\mr e^{-t_ix}\d\Pi_{\beta,w}(x)\right).
\end{align}
Our goal is to make \eqref{Equation: First Term of Inside} and \eqref{Equation: Second Term of Inside} vanish along a 
subsequence. Note that for \eqref{Equation: Second Term of Inside}, using the uniformly bounded variance assumption and the vanishing covariance assumption in 
Lemma \ref{Lemma: New Condition}, the sequence of random variables $\mr{Tr}\left[\mr e^{-t_i H_{\beta,w}}\right]-\mbf{E}\left[\mr{Tr}[\mr e^{-t_iH_{\beta,w}}]\right]$ satisfies precisely the 
assumptions of Lemma \ref{Lemma: WLLN}. Therefore, we can always choose
a sparse enough sequence of $t_i$'s such that
$$\frac{1}{n}\sum_{i=1}^n\left(\mr{Tr}\left[\mr e^{-t_i H_{\beta,w}}\right]-\mbf{E}\left[\mr{Tr}[\mr e^{-t_iH_{\beta,w}}]\right]\right)\rightarrow 0$$
as $n\to \infty$ in probability. It follows that there exists a sequence of positive integers $(n_k)_{k\geq 1}$ so that 
$$\frac{1}{n_k}\sum_{i=1}^{n_k}\left(\mr{Tr}\left[\mr e^{-t_i H_{\beta,w}}\right]-\mbf{E}\left[\mr{Tr}[\mr e^{-t_iH_{\beta,w}}]\right]\right)\rightarrow 0$$
as $k\to \infty$ almost surely. For \eqref{Equation: First Term of Inside}, consider the difference
\begin{align}
  \label{Equation: First Term w/o Average}
  \int_B 1\d\Pi_{\beta,w}(x)-\int_B \mr e^{-t_ix}\d\Pi_{\beta,w}(x)
\end{align}
for a fixed $i$. Since $B$ is bounded above, there exists some $C\in(0,\infty)$ such that $C=\sup(B)$.
Moreover, by Corollary \ref{Corollary: SAO}, we have that 
$\Lambda_{\beta,w}(1)>-\infty$ and $\Pi_{\beta,w}(B)<\infty$ almost surely. Thus on this event,
\eqref{Equation: First Term w/o Average} is bounded above by
\begin{multline}
\label{Equation: Max on the boundary}
\left|\int_B 1\d\Pi_{\beta,w}(x)-\int_B \mr e^{-t_ix}\d\Pi_{\beta,w}(x)\right|
\leq\int_B |1-\mr e^{-t_ix}|\d\Pi_{\beta,w}(x)
\\\leq\Pi_{\be,w}(B)\sup_{x\in B\cap\{\La_{\be,w}(k):k\geq1\}}|1-\mr e^{-t_ix}|
\leq\Pi_{\beta,w}(B)\max\big\{|1-\mr e^{-t_i \Lambda_{\beta,w}(1)}|,|1-\mr e^{- t_i C}|\big\},
\end{multline}
where the last inequality comes from the fact that the maximum of $x\mapsto|1-\mr e^{-t_ix}|$ on any interval
is achieved on the boundary of that interval.
The right-hand side of \eqref{Equation: Max on the boundary}
goes to zero almost surely as $i\rightarrow\infty$. Since $\int_B \mr e^{-t_ix}\d\Pi_{\beta,w}(x)$ converges to $\int_B1\d\Pi_{\beta,w}(x)$ almost surely as $i\to\infty$, so does the average $\frac{1}{n_k}\sum_{i=1}^{n_k}\int_B \mr e^{-t_ix}\d\Pi_{\beta,w}(x)$
as $k\to\infty$. We have now shown that 
\eqref{Equation: First Term of Inside} and \eqref{Equation: Second Term of Inside} both vanish along this subsequence $n_k$, and therefore we have 
\begin{align}
\label{Equation: Inside wrt Outside}
  \Pi_{\beta,w}(B) = \lim_{k\rightarrow\infty}\frac{1}{n_k}\sum_{i=1}^{n_k}\left(\mbf{E}[\mr{Tr}[\mr e^{-t_iH_{\beta,w}}]]-\int_{B^c}\mr e^{-t_ix}\d\Pi_{\beta,w}(x)\right)\qquad\text{almost surely}.
\end{align}
At this point, if we define the deterministic measurable map $F_B:\mc N^\#_{B^c}\to\mbb R$ as
\[F_B(\mu):=\begin{cases}
\displaystyle \lim_{k\rightarrow\infty}\frac{1}{n_k}\sum_{i=1}^{n_k}\left(\mbf{E}[\mr{Tr}[\mr e^{-t_iH_{\beta,w}}]]-\int_{B^c}\mr e^{-t_ix}\d\mu(x)\right)&\text{if the limit exists}\\
0&\text{otherwise}
\end{cases}\]
for all $\mu\in\mc N^\#_{B^c}$, then we get the almost-sure equality \eqref{eq:PiRigidity} by \eqref{Equation: Inside wrt Outside}.

\subsubsection{Proof of Lemma~\ref{Lemma: WLLN}}
\label{Section: WLLN}
Set $r_1 = 1$. Inductively, since $\mbf{E}[X_{r_i}X_j]\rightarrow 0$ as $j\rightarrow\infty$ for each $r_1,\ldots,r_k$, 
 there exists $r_{k+1}>r_k$ so that $\mbf{E}[X_{r_i}X_{r_{k+1}}]\leq 1/(k+1)$ for all $i=1,\ldots,k$. This way, we obtain a subsequence 
 that satisfies $\mbf{E}[X_{r_i}X_{r_j}]\leq 1/i\leq 1/(i-j)$ for every $i>j$. Then, for any fixed $N\in\mbb N$,
  \begin{align*}
    \mbf{E}\left[\left(\frac{X_{r_1}+\cdots +X_{r_n}}{n}\right)^2\right]&=\frac{1}{n^2}\sum_{1\leq i,j\leq n}\mbf{E}[X_{r_i}X_{r_j}]\\
    &\leq\frac1{n^2}\sum_{i=1}^nC+\frac{1}{n^2}\sum_{\substack{1\leq i,j\leq n\\1\leq|i-j|\leq N}}\frac{1}{|i-j|}+\frac{1}{n^2}\sum_{\substack{1\leq i,j\leq n\\|i-j|>N}}\frac1{|i-j|}\\
    &\leq\frac{C}{n}+\frac{2N}{n}+\frac{1}{N}.
  \end{align*}
  The right-hand side of the above inequality converges to $1/N$ as $n\to\infty$, which can be made arbitrarily small by taking $N\to\infty$.
  Therefore, $(X_{r_1}+\cdots +X_{r_n})/n$ converges to zero in $L^2$, hence in probability.

\subsection{Proof of Lemma \ref{Lemma: Covariance Formula}}
\label{Section: Covariance Formula}
This follows from a similar calculation as \cite[Lemma~4.5]{GLGL} using the Feynman-Kac formula of $\mr e^{-tH_{\beta,w}}$ as shown in Theorem~\ref{Theorem: Trace Identity}. We first write 
\[\mbf{Cov}\big[\mr{Tr}[\mr e^{-tH_{\be,w}/2}],\mr{Tr}[\mr e^{-uH_{\be,w}/2}]\big] 
= \mbf{E}[\mr{Tr}[\mr e^{-tH_{\be,w}/2}]\cdot\mr{Tr}[\mr e^{-uH_{\be,w}/2}]]-\mbf{E}[\mr{Tr}[\mr e^{-tH_{\be,w}/2}]]\cdot\mbf{E}[\mr{Tr}[\mr e^{-uH_{\be,w}/2}]].\]
Then, by \eqref{Equation: Trace Identity} and Tonelli's theorem, we have
\begin{align*}
	\mbf{E}[\mr{Tr}[\mr e^{-tH_{\be,w}/2}]]
	&=\int_0^\infty \frac{1+\mr e^{-2x^2/t}}{\sqrt{2\pi t}}\mbf E\left[\mbf E\left[\mr e^{-\langle L_t(X^{x,x}_t),\frac12V\rangle-\frac12\xi_\be(L_t(X^{x,x}_t))-w\mf L_t^0(X^{x,x}_t)}\Big|\xi_\be\right]\right]\d x\\
	&= \int_0^\infty \frac{1+\mr e^{-2x^2/t}}{\sqrt{2\pi t}}\mbf E\left[\mbf E\left[\mr e^{-\langle L_t(X^{x,x}_t),\frac12V\rangle-\frac12\xi_\be(L_t(X^{x,x}_t))-w\mf L_t^0(X^{x,x}_t)}\bigg|X^{x,x}_t\right]\right]\d x\\
	&= \int_0^\infty \frac{1+\mr e^{-2x^2/t}}{\sqrt{2\pi t}}\mbf E\left[\mr e^{-\langle L_t(X^{x,x}_t),\frac12V\rangle+\frac{1}{2\be}\left\Vert L_t(X^{x,x}_t)\right\Vert_2^2-w\mf L_t^0(X^{x,x}_t)}\right]\d x.
\end{align*}
Similarly, 
\[	\mbf{E}[\mr{Tr}[\mr e^{-uH_{\be,w}/2}]]
	= \int_0^\infty \frac{1+\mr e^{-2y^2/t}}{\sqrt{2\pi t}}\mbf E\left[\mr e^{-\langle L_u(\bar{X}^{y,y}_u),\frac12V\rangle+\frac{1}{2\be}\left\Vert L_u(\bar{X}^{y,y}_u)\right\Vert_2^2-w\mf L_u^0(\bar{X}^{y,y}_u)}\right]\d y\]
and
\[	\mbf{E}[\mr{Tr}[\mr e^{-tH_{\be,w}/2}]\cdot \mr{Tr}[\mr e^{-uH_{\be,w}/2}]]
	=\int_{(0,\infty)^2} \mc P_{t,u}(x,y)\mbf E\left[\mr e^{\mc A_{t,u}(x,y)+\mc B_{t,u}(x,y)+\mc C_{t,u}(x,y)+ \mc D_{t,u}(x,y)}\right]\d x\dd y.\]
Subtracting $\mbf{E}[\mr{Tr}[\mr e^{-tH_{\be,w}/2}]]\cdot \mbf{E}[\mr{Tr}[\mr e^{-uH_{\be,w}/2}]]$ from the above expression we complete the proof of Lemma \ref{Lemma: Covariance Formula}.

\begin{notation}
	Throughout the next few subsections we will use $C,C',c,c'$ to denote absolute constants that are independent of $t$. Their exact values may change from line to line.
\end{notation}
\subsection{Proof of Lemma \ref{Lemma: Self-Intersection}}
\label{Section: Self-Intersection}
The proof of this Lemma shares similar ideas as in \cite[Lemma~4.6]{GLGL}. To avoid repetitions, we will omit details in few occasions and refer to relevant places in the proof of \cite[Lemma~4.6]{GLGL}.
We start with the bound for the $\mc C_{t,u}(x,y)$ part. By independence we have 
\[
\mbf E\left[\mr e^{4\mc C_{t,u}(x,y)}\right] = \mbf E\left[\mr e^{\frac{2}{\beta}\Vert L_t(X_t^{x,x})\Vert_2^2}\right]\cdot \mbf E\left[\mr e^{\frac{2}{\beta}\Vert L_u(\bar{X}_u^{y,y})\Vert_2^2}\right].
\]
We will show that there exists constants $C',c'>0$ (that only depend on $\be$) such that for all $t>0$,
\[
	\sup_{x>0}\mbf E\left[\mr e^{\frac{2}{\beta}\Vert L_t(X_t^{x,x})\Vert_2^2}\right]\leq C \mr e^{c'\left(t^{3/2}+t^3\right)}.
\]
This is a slightly stronger statement comparing to the one in \cite[Lemma~4.6]{GLGL} which only establishes the bound for $t>0$ sufficiently small. However the original proof still works in this new set up. We first condition on the value $X_{t}^{x,x}(t/2)$ and use Doob's $h$-transform to write
\[\mbf E\left[\mr e^{\frac{2}{\beta}\Vert L_t(X_t^{x,x})\Vert_2^2}\right]= \int_0^{\infty} \mbf E\left[\mr e^{\frac{2}{\beta}\Vert L_t(X_t^{x,x})\Vert_2^2}\big| X_t^{x,x}(t/2)=y\right]\frac{\Pi_{X}(t/2;x,y)\Pi_X(t/2;y,x)}{\Pi_{X}(t;x,x)}\dd y,\]
where 
$\Pi_{X}(t;x,y) = \frac{\mr e^{-(x-y)^2/2}+\mr e^{-(x+y)^2/2}}{\sqrt{2\pi t}}$ is the transition kernel for the reflected Brownian motion. Now by \cite[(4.22)]{GLGL} we have
\[
	\mbf E\left[\mr e^{\frac{2}{\beta}\Vert L_t(X_t^{x,x})\Vert_2^2}\big| X_t^{x,x}(t/2)=y\right]\leq \mbf E\left[\mr e^{\frac{8}{\beta}\Vert L_{t/2}(X_{t/2}^{x,y})\Vert_2^2}\right].	
\]
On the other hand we have the elementary bounds 
\[
	\Pi_{X}(t;x,x)\geq \frac{1}{\sqrt{2\pi t}},\quad \Pi_{X}(t/2;y,x)\leq \frac{2}{\sqrt{\pi t}},
\]
for all $t,x,y>0$. Hence 
\[
 \sup_{t>0}\sup_{x,y>0} \frac{\Pi_X(t/2;y,x)}{\Pi_{X}(t;x,x)}\leq 2\sqrt{2}.
\]
Combine the bounds together we see that 
\begin{equation}
	\begin{aligned}
	\mbf E\left[\mr e^{\frac{2}{\beta}\Vert L_t(X_t^{x,x})\Vert_2^2}\right]&\leq 2\sqrt{2}\int_0^{\infty}\mbf E\left[\mr e^{\frac{8}{\beta}\Vert L_{t/2}(X_{t/2}^{x,y})\Vert_2^2}\right]\Pi_{X}(t/2;x,y)\d y\\
	&= 2\sqrt{2}\mbf E\left[\mr e^{\frac{8}{\beta}\Vert L_{t/2}(X^x)\Vert_2^2}\right].
	\end{aligned}
\end{equation}
Hence it suffices to show that for all $t>0$ we have 
\[
	\sup_{x>0} \mbf E\left[\mr e^{\frac{8}{\beta}\Vert L_{t/2}(X^x)\Vert_2^2}\right] \leq C \mr e^{c'\left(t^{3/2}+t^3\right)},
\]
where $C',c'>0$ only depend on $\be$.
For this we couple $X$ with a standard Brownian motion $B$ so that $X^x(t)=|B^x(t)|$. Then for every $a>0$ we have $L_t^a(X^x)=L_t^a(B^x)+L_t^{-a}(B^x)$. Hence
\begin{align}
\label{Equation: Scaling In Distribution}
	\Vert L_t(X^x)\Vert_2^2=\int_{0}^\infty L_t^a(X^x)^2\d a\leq 2\int_\mbb R L_t^a(B^x)^2\d a
	\deq2\int_\mbb R L_t^a(B^0)^2\d a\deq 2t^{3/2}\Vert L_1(B^0)\Vert_2^2,
\end{align}
where the second-to-last equality in distribution follows from the change of variables $a\mapsto a-x$,
and the last equality in distribution follows from the Brownian scaling identity (e.g., \cite[Propositions 2.3.3 and 2.3.5]{Chen}).
By \cite[Theorem 4.2.1]{Chen} we have the following sub-Gaussian tail bound for $\Vert L_1(B^0)\Vert_2^2$: 
\[
	\mbf P[\Vert L_1(B^0)\Vert_2^2>s]\leq C\mr e^{-cs^2}, 
\]
for absolute constants $c,C>0$.
Recall that a real random variable $Z$ is sub-Gaussian if and only if there 
exist constants $B,b>0$ such that
\begin{align}
\label{Equation: Sub-G}
\mbf E[\mr e^{t(Z-\mbf EZ)}]\leq B\mr e^{bt^2}\qquad\text{for all }t\in\mbb R
\end{align}
(e.g. \cite[Proposition 2.5.2-(v)]{Vershynin}).
Thus, applying this property and \eqref{Equation: Scaling In Distribution} we conclude that
\begin{multline*}
\sup_{x>0} \mbf E\left[\mr e^{\frac{8}{\beta}\Vert L_{t/2}(X^x)\Vert_2^2}\right]
\leq\mbf E\left[\mr e^{\frac{4 \sqrt{2} t^{3/2}}{\beta}\Vert L_{1}(B^0)\Vert_2^2}\right]\\
=\mbf E\left[\mr e^{\frac{4 \sqrt{2} t^{3/2}}{\beta}\mbf E\Vert L_{1}(B^0)\Vert_2^2+\frac{4 \sqrt{2} t^{3/2}}{\beta}(\Vert L_{1}(B^0)\Vert_2^2-\mbf E\Vert L_{1}(B^0)\Vert_2^2)}\right]
\leq C'\mr e^{c'(t^{3/2}+t^3)}
\end{multline*}
for constants $c',C'>0$ that only depend on $\be$. Hence 
\[
	\limsup_{t\to0}\sup_{x,y>0}\mbf E\left[\mr e^{4\mc C_{t,u}(x,y)}\right]\leq c<\infty.
\]
The proof for the $\mc B_{t,u}(x,y)$ part is similar, see also \cite[Lemma 5.6]{GL0} for a similar statement.

\subsection{Proof of Lemma \ref{Lemma: Gaussian Term}}
\label{Section: Gaussian Term}
By the definition of $\mc D_{t,u}$ (see Section~\ref{Section: Semigroup Approach}), we have $\mc D_{t,u}\neq 0$ if and only if $\langle L_t(X_t^{x,x}),L_u(\bar{X}_{u}^{y,y})\rangle\neq 0$ where $\langle \cdot,\cdot\rangle$ is the standard $L^2$ inner product. Hence by Cauchy-Schwarz inequality we have
\begin{equation}
	\label{Equation: Inserting indicator}
\begin{aligned}
	\mbf E\left[\left(\mr e^{\mc D_{t,u}(x,y)}-1\right)^4\right]^{1/4}
	&=\mbf E\left[\left(\mathbf{1}_{\{\langle L_t(X_t^{x,x}),L_u(\bar{X}_{u}^{y,y})\rangle\neq 0\}}\left(\mr e^{\mc D_{t,u}(x,y)}-1\right)\right)^4\right]^{1/4}\\
	&\leq \mbf E\left[\left(\mr e^{\mc D_{t,u}(x,y)}-1\right)^8\right]^{1/8}\cdot \mbf P\left[\langle L_t(X_t^{x,x}),L_u(\bar{X}_{u}^{y,y})\rangle\neq 0\right]^{1/8}.
\end{aligned}
\end{equation}
To estimate $\mbf P\left[\langle L_t(X_t^{x,x}),L_u(\bar{X}_{u}^{y,y})\rangle\neq 0\right]$ we first couple the processes $X$ and $\bar{X}$ with two independent standard Brownian motions $B$ and $\bar{B}$ so that 
\begin{equation}
	\label{Equation: Coupling}
  X^{x}(s)= |x+B^0(s)|,\quad \bar{X}^{y}(s)= |y+\bar{B}^0(s)|.
\end{equation}
Then by \cite[(5.9)]{GL0} we have for any nonnegative path functional $F$,
\begin{equation}
	\label{Equation: Coupling bound}
	\mbf E[F(X_t^{x,x})]\leq 2\mbf E[F(|x+B_t^{0,0}|)].
\end{equation}
In particular
\[
\mbf P\left[\langle L_t(X_t^{x,x}),L_u(\bar{X}_{u}^{y,y})\rangle\neq 0\right]\leq 2\mbf P\left[\langle L_t(|x+B_t^{0,0}|),L_u(|y+\bar{B}_{u}^{0,0}|)\rangle\neq 0\right].
\]
Now it is easy to check the following inclusion of events:
\begin{equation}
	\label{Equation: union bound}
	\begin{aligned}
	&\{\langle L_t(|x+B_t^{0,0}|),L_u(|y+\bar{B}_{u}^{0,0}|)\rangle\neq 0\}\\
	&\subset \{\text{the range of $(|x+B_t^{0,0}(s)|)_{s\in [0,t]}$ and $(|y+\bar{B}_u^{0,0}(s)|)_{s\in [0,u]}$ intersect}\}\\
	&\subset \{\max_{s\in [0,t]}|B_t^{0,0}(s)|\geq \frac{|x-y|}{2}\}\cup\{\max_{s\in [0,u]}|\bar{B}_u^{0,0}(s)|\geq \frac{|x-y|}{2}\}.
	\end{aligned}
\end{equation}
By Brownian scaling and the fact that maximum of Brownian bridge have sub-Gaussian tails we have 
\[
   \mbf P[\max_{s\in [0,t]}|B_t^{0,0}(s)|\geq \frac{|x-y|}{2}]= \mbf P[\max_{s\in [0,1]}|B_1^{0,0}(s)|\geq \frac{|x-y|}{2t^{1/2}}]\leq C\mr e^{-c(x-y)^2/t}
\]
for some $C,c>0$ independent of $t$.
Therefore,
\begin{align*}
	\mbf P\left[\langle L_t(|x+B_t^{0,0}|),L_u(|y+B_{u}^{0,0}|)\rangle\neq 0\right]
	&\leq C\left(\mr e^{-c(x-y)^2/t}+\mr e^{-c(x-y)^2/u}\right)\\
	&\leq C'\mr e^{-c'(x-y)^2}
\end{align*}
for any $t\in (0,1]$, where $c'=c'(u)= c\cdot\min\{1,1/u\}$. Combining this with \eqref{Equation: Inserting indicator} we complete the proof of Lemma \ref{Lemma: Gaussian Term}.

\subsection{Proof of Lemma \ref{Lemma: Potential Term}}
\label{Section: Potential Term}
By the definition of $\mc A_{t,u}(x,y)$ and basic properties of local time we have 
\begin{align*}
	\mc A_{t,u}(x,y) = -\langle L_t(X_t^{x,x})+L_u(\bar{X}_u^{y,y}),\tfrac{1}{2}V\rangle= -\frac{1}{2}\left(\int_{0}^t X_t^{x,x}(s)\d s+\int_{0}^u \bar{X}_u^{y,y}(s)\d s\right).
\end{align*}
Under the same coupling as in \eqref{Equation: Coupling} we have by \eqref{Equation: Coupling bound}
\[
    \mbf E\left[\mr e^{-2\int_{0}^t X_t^{x,x}(s)\d s}\right]\leq 2\mbf E\left[\mr e^{-2\int_{0}^t |x+B_t^{0,0}(s)|\d s}\right]\leq 2 \mr e^{-2xt}\mbf E\left[\mr e^{2\int_0^t|B_t^{0,0}(s)|\d s}\right].
\]
Now by a change of variable $s\mapsto st$ and Brownian scaling we have 
\begin{align*}
	\mbf E\left[\mr e^{2\int_0^t|B_t^{0,0}(s)|\d s}\right]= \mbf E\left[\mr e^{2t^{3/2}\int_0^1|B_1^{0,0}(s)|\d s}\right]\leq \mbf E\left[\mr e^{2t^{3/2}\max_{s\in [0,1]}|B_1^{0,0}(s)|}\right].
\end{align*}
Since the maximum of the Bessel bridge $|B_1^{0,0}(s)|$ has a sub-Gaussian tail,
we have by \eqref{Equation: Sub-G} that 
\begin{align*}
	\mbf E\left[\mr e^{2t^{3/2}\max_{s\in [0,1]}|B_1^{0,0}(s)|}\right]\leq  C\mr e^{c(t^{3/2}+t^3)}
\end{align*}
for some constants $c,C>0$ independent of $t$. Thus by independence we have 
\begin{align*}
	\mbf E \left[\mr e^{4\mc A_{t,u}(x,y)}\right]^{1/4}
	&= \left(\mbf E \left[\mr e^{-2\int_{0}^t X_t^{x,x}(s)\d s}\right]\cdot \mbf E \left[\mr e^{-2\int_{0}^u \bar{X}_u^{y,y}(s)\d s}\right]\right)^{1/4}\\
	&\leq C \mr e^{-\frac{1}{2}xt-\frac{1}{2}yu+c(t^{3/2}+t^3+u^{3/2}+u^3)}
\end{align*}
Therefore, we can find constants $C,c>0$ such that for any $t\in (0,1]$ and $u>0$ fixed we have 
\begin{align*}
	\int_{(0,\infty)^2}\mbf E\left[\mr e^{4\mc A_{t,u}(x,y)}\right]^{1/4}\mr e^{-(x-y)^2/c}\d x\dd y
	&\leq C\int_{(0,\infty)^2}\mr e^{-yu/2-(x-y)^2/c}\d x\dd y<\infty.
\end{align*}
This completes the proof of Lemma \ref{Lemma: Potential Term}.

\subsection{Proof of Lemma \ref{Lemma: Vanishing Term}}
\label{Section: Vanishing Term}
This is again a similar calculation as \cite[Lemma 4.8]{GLGL}. Using Cauchy-Schwarz and the elementary inequality $|\mr e^z -1|\leq |z|\mr e^{|z|}$ we have 
\begin{equation}
	\label{Equation: Estimate D}
	\left|\mr e^{\mc D_{t,u}(x,y)}-1\right|\leq \frac{1}{\beta}\Vert L_t(X_t^{x,x})\Vert_2\Vert L_u(\bar{X}_u^{y,y})\Vert_2 \cdot \mr e^{\frac{1}{\beta}\Vert L_t(X_t^{x,x})\Vert_2\Vert L_u(\bar{X}_u^{y,y})\Vert_2}.
\end{equation}
By Cauchy-Schwarz inequality we have 
\[
	\mbf E\left[\Vert L_t(X_t^{x,x})\Vert_2^8\mr e^{\frac{8}{\beta}\Vert L_t(X_t^{x,x})\Vert_2}\right]\leq \left(\mbf E\Vert L_t(X_t^{x,x})\Vert_2^{16}\right)^{1/2}\left(\mbf E\mr e^{\frac{16}{\beta}\Vert L_t(X_t^{x,x})\Vert_2}\right)^{1/2}.
\]
Now a similar argument as in Lemma \ref{Lemma: Self-Intersection} shows that 
\[
	\mbf E\Vert L_t(X_t^{x,x})\Vert_2^{16}\leq C\mbf E\Vert L_{t/2}(B^x)\Vert_2^{16}= C \left(\frac{t}{2}\right)^{12}\mbf E\Vert L_1(B^0)\Vert_2^{16}\leq C' t^{12},
\]
$B^0$ is a standard Brownian motion and we couple $X$ and $B$ as in \eqref{Equation: Coupling}. Similarly
\begin{align*}
	\mbf E\mr e^{\frac{16}{\beta}\Vert L_t(X_t^{x,x})\Vert_2}\leq C\mbf E\mr e^{\frac{64}{\beta}\Vert L_{t/2}(X^x)\Vert_2}\leq C \mbf E\mr e^{c\Vert L_{t/2}(B^x)\Vert_2}=C' \mbf E\mr e^{c't^{3/4}\Vert L_{1}(B^0)\Vert_2}.
\end{align*}
By \cite[Theorem 4.2.1]{Chen} we have
\[
	\mbf P[\Vert L_1(B^0)\Vert_2>s]\leq C\mr e^{-cs^4}.
\]
Hence using once again the sub-Gaussian bound \eqref{Equation: Sub-G} we get
\[
\mbf E\mr e^{\frac{16}{\beta}\Vert L_t(X_t^{x,x})\Vert_2}\leq C \mbf E\mr e^{\frac{c}{\beta}t^{3/4}\Vert L_{1}(B^0)\Vert_2}\leq C'\mr e^{c'\left(t^{3/2}+t^3\right)}.
\]
Combine the estimates together we see that 
\begin{equation}
	\label{Equation: Estimate D II}
	\mbf E\left[\Vert L_t(X_t^{x,x})\Vert_2^8\mr e^{\frac{8}{\beta}\Vert L_t(X_t^{x,x})\Vert_2}\right]\leq C t^{6}\mr e^{c\left(t^{3/2}+t^3\right)}.
\end{equation}
Finally using \eqref{Equation: Estimate D}, \eqref{Equation: Estimate D II} and independence we conclude that 
\[
	\left(\mbf E\left|\mr e^{\mc D_{t,u}(x,y)}-1\right|^8\right)^{1/8}\leq C_{u,\beta}\cdot t^{3/4},
\]
uniformly for all $t\in (0,1]$ where the constant $C_{u,\beta}$ depends only on $u,\beta>0$ but not $t$. This completes the proof of Lemma \ref{Lemma: Vanishing Term}.

\subsection*{Acknowledgments}

Y. Liao is supported by EPSRC grant EP/R024456/1. Part of this research was conducted while P. Y. Gaudreau Lamarre and W. Li were involved in
the University of Chicago Mathematics REU 2021. The organizers of the REU are gratefully
acknowledged for fostering a productive research environment. P. Ghosal wishes to thank Amol Aggarwal, Benjamin Mckenna and Jos\'e Ram\'irez for many useful conversations.

\end{document}